\documentclass[11pt]{article}

\usepackage{amsmath}
\usepackage{amsthm}
\usepackage{amssymb}
\usepackage[margin=1in]{geometry}
\usepackage{url}

\theoremstyle{plain}
\newtheorem{theorem}{Theorem}[section]

\newtheorem{corollary}[theorem]{Corollary}
\theoremstyle{definition}

\theoremstyle{remark}

\newcommand{\Alb}{\operatorname{Alb}}
\newcommand{\rank}{\operatorname{rank}}

\title{New excluded minors for the class $\mathcal{M}_3$ of regular matroids}
\author{Booker Smith\\
\texttt{bookers7@gmail.com}}
\date{19 July 2026}

\begin{document}
\maketitle

\begin{abstract}
Engel, de Gaay Fortman, and Schreieder attach to each prime $\ell$ a minor-closed
class $\mathcal{M}_\ell$ of regular matroids, whose excluded minors govern the
failure of the integral Hodge conjecture for curve classes on very general
principally polarized abelian varieties. For $\ell = 2$ the class is exactly the
cographic matroids, with excluded minors $M(K_5)$ and $M(K_{3,3})$ by Tutte's
theorem; for the odd prime $\ell = 3$ the only excluded minor explicitly
identified so far is
$M(K_{3,5})$ [EGFS, Prop.~8.11], and the general characterization is
[EGFS, Problem~8.8]. We exhibit five new excluded minors for $\mathcal{M}_3$, of
ranks $8$, $9$, $9$, $9$, and $10$. None contains $M(K_{3,5})$ or any of the
other four as a minor, so the excluded-minor list of Problem~8.8 has at least six
members, with excluded minors at every rank from $7$ through $10$. Combined with
the non-membership $M(K_7) \notin \mathcal{M}_3$ established by the same authors
[EGFS2] and the minor-closedness of $\mathcal{M}_3$, some minor of $M(K_7)$ of
rank between $4$ and $6$ is a further excluded minor, distinct from these six;
the list therefore has at least seven members, six of them explicitly
identified. The list is
structurally diverse: the rank-$8$ and rank-$10$ minors are graphic matroids of
bipartite apex constructions with automorphism group $S_4$ (the $3$-cube plus a
one-sided apex, and the subdivided $K_4$ plus an apex), whereas the three
rank-$9$ minors are graphic matroids of non-bipartite, non-planar graphs with
automorphism groups of order at most $4$ and no such apex structure. Every
verdict is certified by a finite $\mathbb{F}_3$-linear-algebra computation with
explicit machine-checkable witnesses, cross-checked by independent
implementations. Combined with a recent result of Engel and Schreieder
[ES, Cor.~2.10], the certificates further show that the five matroids lie outside
the larger class $\widetilde{\mathcal{M}}_3$ of [EGFS, Rmk.~8.5]. We further
record what are, to our knowledge, the first $\ell = 5$ data on this corank-$8$
slice: $M(K_{3,5})$ and the rank-$8$ and rank-$10$ minors lie in $\mathcal{M}_5$,
each by an explicit certificate, so none of these three contributes a factor of
$5$ to its radical distance; combined with the rank bound of [EGFS, Rmk.~8.6],
the $\ell = 5$ certificate determines the radical distance
$d(M(K_{3,5})) = 6$ exactly. The $\ell = 5$ status of the three rank-$9$ minors
has not yet been computed. This complements the recent theorem of Engel and
Schreieder [ES] that the minimal algebraic multiple of the minimal class on a
very general principally polarized abelian variety is divisible by all primes
$p \le (k+1)/2$, which resolves [EGFS, Problem~8.13].
\end{abstract}

\section{Setup}\label{sec:setup}

We work in the framework of [EGFS]. To a regular matroid $(R,S)$ of rank $g$ on
ground set $S$, with a totally unimodular integral realization
$S \hookrightarrow U^{*}$ (equivalently a matrix
$A = (I_g \mid D) \in \mathbb{Z}^{g \times |S|}$, with row space
$U = U(R) \subset \mathbb{Z}^{S}$; for a graphic matroid $M(G)$, $U$ is the cut
lattice of $G$), [EGFS] associate the
following combinatorial data.

\medskip

\noindent\textbf{The Albanese graph} [EGFS, Def.~5.3]. For non-negative integers
$j \le r$, the $(\ell^{r}, \ell^{j})$-Albanese graph
$\Alb_{\ell^{r},\ell^{j}}(R)$ is the oriented $S$-colored graph with vertex set
\[
V \;=\; \mathbb{Z}^{S}/(\ell^{j} U + \ell^{r}\mathbb{Z}^{S}),
\]
in which $[v]$ and $[w]$ are joined by an oriented edge of color $s \in S$,
directed from $[v]$ to $[w]$, precisely when $[w] = [v + e_s]$. We write
$\Alb_{\ell}(R) := \Alb_{\ell^{1},\ell^{0}}(R)$, whose vertex set is
$\mathbb{F}_\ell^{S}/\bar U$ (with $\bar U = U \bmod \ell$, of dimension $g$), so
that $|V| = \ell^{\,|S| - g}$ and each vertex has one out-edge of each color.

\medskip

\noindent\textbf{The color profile} [EGFS, Def.~5.1] is the linear map
$\lambda \colon C_1(\Alb_\ell(R),\Lambda) \to \Lambda^{S}$ sending an oriented
edge of color $s$ to the basis vector $e_s$; for a $1$-chain $b$ we write
$\lambda_s(b)$ for its $s$-component.

\medskip

\noindent\textbf{Solutions and membership} [EGFS, Def.~5.8, Def.~8.3]. A
$\Lambda$-solution of $(R,S)$ in a graph $G$ is a family $(b_s)_{s\in S}$ of
$1$-chains, $b_s$ supported on the color-$s$ edges, such that $\sum_{s} c_s\, b_s$
is closed whenever $\sum_s c_s e_s \in U_\Lambda$. The matroid $R$ lies in the
class
\[
\mathcal{M}_\ell \;:=\; \bigl\{\, R \text{ regular} \;:\; R \text{ admits a }
\mathbb{Z}/\ell\text{-solution in } \Alb_{\ell}(R) \text{ with constant nonzero
color profile} \,\bigr\},
\]
i.e.\ a solution with $\lambda_s(b_s) = c$ for a single $c \in (\mathbb{Z}/\ell)^{\times}$
and all $s$.

By [EGFS, Prop.~7.2] each $\mathcal{M}_\ell$ is closed under taking minors;
within graphic matroids, the Robertson--Seymour graph-minor theorem therefore
guarantees a characterization by finitely many excluded graph minors. For the
full class of regular matroids finiteness is expected but open ---
$\mathcal{M}_\ell$ ``is expected to be determined by a finite list of excluded
minors'' via the conjectural extension of the Robertson--Seymour theorem to
regular matroids [EGFS, \S8.5] --- and Problem~8.8 below asks for the
characterization. Two facts frame the present note.

\begin{itemize}
\item \textbf{$\ell = 2$ is classical.} $\mathcal{M}_2$ is the class of cographic
matroids [EGFS, Thm.~7.1]; a graphic matroid $M(G)$ is cographic iff $G$ is
planar. Tutte's theorem gives the two excluded minors $M(K_5)$ and $M(K_{3,3})$.

\item \textbf{$\ell = 3$ had, until now, a single explicitly identified excluded
minor.}
$M(K_{3,5})$ (rank $7$) satisfies $M(K_{3,5}) \notin \mathcal{M}_3$ while every
single-element minor lies in $\mathcal{M}_3$, hence is an excluded minor for
$\mathcal{M}_3$ [EGFS, Prop.~8.11]. As of both [EGFS] and the subsequent work
[ES] of Engel and Schreieder, $M(K_{3,5})$ remains the only excluded minor
explicitly identified for $\mathcal{M}_3$ in the literature: [ES] proves the
further non-memberships $M(K_{2p+1}) \notin \mathcal{M}_p$ and
$M(K_{3,2p-1}) \notin \mathcal{M}_p$ for every prime $p$ (its Theorems~1.5
and~1.7; the case $p = 3$ recovers $M(K_7) \notin \mathcal{M}_3$ --- first
established in the follow-up [EGFS2, Cor.~3.9] --- and
$M(K_{3,5}) \notin \mathcal{M}_3$), but it identifies no excluded minors and
leaves Problem~8.8 open. The corresponding open problem is:
\end{itemize}

\begin{quote}
\textbf{Problem 8.8 [EGFS].} \emph{Let $\ell$ be an odd prime. Characterize the
subclass $\mathcal{M}_\ell$ of the class of regular matroids via a finite list of
excluded minors.}
\end{quote}

The radical distance
$d(R) := \operatorname{lcm}\{\ell \text{ prime} : R \notin \mathcal{M}_\ell\}$
[EGFS, Def.~8.9] measures how far $R$ is from cographic; it is finite and divides
the product of the primes below $\rank(R)$, and $d(R) = 1$ iff $R$ is cographic.
Excluded minors for $\mathcal{M}_\ell$ have rank $> \ell$
[EGFS, Rmk.~8.6 / Thm.~1.10].

\section{The results}\label{sec:result}

We describe five graphs $G_1, \dots, G_5$; all are non-planar of girth $4$ with
vertex connectivity exactly $3$, and each graphic matroid $M(G_i)$ is regular,
simple, and of corank $8$.

\medskip

\noindent\textbf{The rank-$8$ graph $G_1$.} Let $G_1$
be the bipartite graph on parts $X = \{1,2,3,4\}$ and $Y = \{5,6,7,8,9\}$ with
edge set
\[
E(G_1) = \{\,15,25,35,\ \ 16,26,46,\ \ 17,37,47,\ \ 28,38,48,\ \ 19,29,39,49\,\}.
\]
Equivalently, $G_1 = K_{4,5}$ with the $4$-edge matching $\{18, 27, 36, 45\}$
deleted (each vertex of $X$ loses one edge, to four of the five vertices of $Y$;
vertex $9$ remains complete to $X$). It has $9$ vertices, $16$ edges, degree
sequence $(4^{5}, 3^{4})$, girth $4$, is $3$-connected and non-planar, with
automorphism group of order $24$, isomorphic to $S_4$. Its graph6 code is
\verb|H?zTbbo| (canonical form \verb|H?Bmtpw|). There is a second useful
description: deleting the unique degree-$4$ vertex of $Y$, namely $9$, leaves the
$3$-regular bipartite graph on $\{1,2,3,4\} \cup \{5,6,7,8\}$, which is $K_{4,4}$
minus a perfect matching, i.e.\ the $3$-cube $Q_3$. Thus
\[
G_1 \;=\; Q_3 \text{ together with an apex joined to all four vertices of one
side of } Q_3 .
\]
$M(G_1)$ has rank $8$ and size $16$.

\medskip

\noindent\textbf{The rank-$10$ graph $G_2$.} Let $G_2$
be the graph obtained from $K_4$ by subdividing each of its six edges once and
joining a new apex vertex to the six subdivision vertices: $11$ vertices, $18$
edges, bipartite with parts of sizes $6$ and $5$, degree sequence
$(6^{1}, 3^{10})$, girth $4$, non-planar, automorphism group $S_4$ (order $24$).
Its graph6 code is \verb|J??FCpSJFw?| (canonical form
\verb|J?hcab?K?^_|). $M(G_2)$ has rank $10$ and size $18$.

\medskip

\noindent\textbf{The three rank-$9$ graphs $G_3, G_4, G_5$.}
Each is a graph on $10$ vertices with $17$ edges, degree sequence
$(4^{4}, 3^{6})$, girth $4$, \emph{non-bipartite} and non-planar; each graphic
matroid has rank $9$ and size $17$. On vertex set $\{0, \dots, 9\}$:
\begin{align*}
E(G_3) &= \{04, 06, 07,\ 15, 16, 17,\ 26, 27, 28,\ 37, 38, 39,\ 48, 49,\ 58, 59,\ 69\},\\
E(G_4) &= \{05, 07, 08,\ 15, 17, 18,\ 26, 27, 29,\ 36, 38, 39,\ 47, 48, 49,\ 56, 59\},\\
E(G_5) &= \{03, 05, 06,\ 14, 16, 18,\ 25, 26, 27,\ 37, 38, 39,\ 47, 49,\ 58, 59,\ 69\}.
\end{align*}
Their graph6 codes are \verb|I?`FF_{F_| (canonical form, computed with
\texttt{bliss}, \verb|I@?LQngt?|),
\verb|I?B@nRWN?| (canonical \verb|I??[rNoy?|), and \verb|ICQf@pSF_|
(canonical \verb|IAGWMUqw_|); the three graphs are pairwise non-isomorphic.
Their automorphism groups are small: $C_2 \times C_2$, $C_2 \times C_2$, and
$C_2$ respectively. Unlike $G_1$ and $G_2$, the three rank-$9$ graphs do not
present themselves as a small structured base plus an apex vertex; they were
found by unstructured census search rather than by construction.

\begin{theorem}\label{thm:main}
Each of the five graphic matroids $M(G_1), \dots, M(G_5)$ is an excluded minor
for the class $\mathcal{M}_3$, and no one of the six matroids
$M(K_{3,5}), M(G_1), \dots, M(G_5)$ contains another as a minor. Moreover each
$M(G_i)$ lies outside the larger class $\widetilde{\mathcal{M}}_3$ of
[EGFS, Rmk.~8.5], and is an excluded minor for $\widetilde{\mathcal{M}}_3$ as
well.
\end{theorem}

The non-memberships, the minimality, and the pairwise independence are
established by explicit certificates in \S\ref{sec:method}, where the
$\widetilde{\mathcal{M}}_3$-strengthening is also proved.

\begin{corollary}\label{cor:six}
The class $\mathcal{M}_3$ is not characterized among regular matroids by the
exclusion of $M(K_{3,5})$ alone. Any finite excluded-minor list for
$\mathcal{M}_3$ (Problem~8.8) contains at least six members: $M(K_{3,5})$ and
$M(G_1), \dots, M(G_5)$.
\end{corollary}

\begin{corollary}\label{cor:rank}
Excluded minors for $\mathcal{M}_3$ occur at every rank from $7$ to $10$:
$M(K_{3,5})$ has rank $7$, and the five matroids above have ranks
$8, 9, 9, 9, 10$.
\end{corollary}

\begin{corollary}\label{cor:seven}
The excluded-minor list of Problem~8.8 for $\mathcal{M}_3$ has at least seven
members, of which the six matroids of Corollary~\ref{cor:six} are the ones
explicitly identified.
\end{corollary}

\noindent Indeed, $M(K_7) \notin \mathcal{M}_3$ [EGFS2, Cor.~3.9]. Since
$\mathcal{M}_3$ is minor-closed [EGFS, Prop.~7.2], some minor of $M(K_7)$ ---
necessarily graphic --- is minor-minimal outside $\mathcal{M}_3$, i.e.\ is an
excluded minor for $\mathcal{M}_3$. Its rank is at most
$\rank M(K_7) = 6$, and at least $4$, since $R \in \mathcal{M}_3$ for every
regular matroid $R$ of rank at most $3$ [EGFS, Rmk.~8.6]. The six excluded
minors of Corollary~\ref{cor:six} have ranks $7$ through $10$, so this seventh
one is distinct from all of them. Computations subsequent to the search of
Section~5 single out a candidate: the graphic matroid of the complete
multipartite graph $K_{1,2,2,2}$ --- equivalently, $K_7$ minus a perfect
matching on six of its seven vertices --- of rank $6$, with $18$ edges and
automorphism group of order $48$. Its nonmembership in $\mathcal{M}_3$ is
certified by an explicit dual witness included in the certificate bundle.
Its $36$ single-element minors have not yet received primal certificates to
the standard of Section~3, so we record $M(K_{1,2,2,2})$ as the candidate
seventh excluded minor: its minor-minimality is supported by the same
computations but is not certified here.

\noindent\textbf{Structural diversity.} The list of explicitly identified
excluded minors is not the work of a single construction. $G_1$ and $G_2$
are bipartite base-plus-apex graphs with automorphism group $S_4$; the three
rank-$9$ graphs are non-bipartite, admit no such decomposition, and have
automorphism groups of order at most $4$. In the systematic search that produced
these results (\S\ref{sec:remarks}), $214$ constructed base-plus-apex candidates
over $24$ distinct bases yielded \emph{zero} non-members of $\mathcal{M}_3$,
while the undirected census of generic corank-$8$ graphs produced the three
rank-$9$ excluded minors. On the present data, a description of the excluded-minor
list of Problem~8.8 by a single structural family is not supported.

\medskip

\noindent\textbf{Consequences for radical distance.} These witnesses give no
new lower bound on the radical-distance function. Each $G_i$ is non-planar, so
each $M(G_i)$ is non-cographic and $M(G_i) \notin \mathcal{M}_2$; together with
$M(G_i) \notin \mathcal{M}_3$ this yields $6 \mid d\bigl(M(G_i)\bigr)$, and hence
$6 \mid d(g)$ for $g = 8, 9, 10$, where
$d(g) := \operatorname{lcm}\{d(R) : R \text{ regular of rank } g\}$. But
$6 \mid d(g)$ is already known for all $g \ge 7$ from
$M(K_{3,5}) \preceq M(K_{3,n})$ (rank $n+2$). Thus the five matroids are
independent witnesses for known divisibilities, and contribute no new numerical
floor. A genuinely new floor from these matroids would require one of them to
lie outside $\mathcal{M}_p$ for some prime $p \ge 5$; by the rank bound of
[EGFS, Rmk.~8.6] ($R \in \mathcal{M}_\ell$ whenever $\ell \ge \rank R$), only
$p = 5$ and $p = 7$ are in question at ranks $8$--$10$. For $M(G_1)$ and
$M(G_2)$ the $p = 5$ possibility has since been tested and fails: both lie in
$\mathcal{M}_5$ (\S\ref{sec:ell5}, by explicit certificates). This leaves
$p = 7$, which we have not attempted, as the only undetermined prime for these
two, so $d(M(G_1))$ and $d(M(G_2))$ each lie in $\{6, 42\}$. For the three
rank-$9$ matroids neither the $\ell = 5$ nor the $\ell = 7$ computation has been
carried out, and we make no claim about them.

\section{The certificates and the method}\label{sec:method}

\noindent\textbf{A linear-algebra formulation of membership.} Scaling the
constant profile to $c = 1$, [EGFS, Def.~8.3] reads: $R \in \mathcal{M}_\ell$ iff
the single inhomogeneous $\mathbb{F}_\ell$-linear system
\[
\text{(closedness)} \quad \textstyle\sum_{s} u_s\, \partial b_s = 0 \ \ \text{for a
basis } \{u\} \text{ of } \bar U, \qquad
\text{(profile)} \quad \lambda_s(b_s) = 1 \ \ \text{for all } s\in S,
\]
in the unknown $1$-chains $(b_s)_{s\in S}$ is solvable. Solvability is decided by
rank comparison during the search; the certificates themselves use no rank
computation. Membership is certified by an explicit solution (a \emph{primal
witness}), and non-membership by an explicit \emph{dual certificate} --- a left-kernel witness
$y$ with $y^{\mathsf T} M = 0$ and $y^{\mathsf T} b \ne 0$ (a Farkas-type
alternative over $\mathbb{F}_3$) --- each checkable by direct matrix--vector
arithmetic over $\mathbb{F}_3$.

\medskip

\noindent\textbf{Non-membership.} Write $C$ for the closedness matrix
($g \cdot 3^{8}$ rows, one per basis element of $\bar U$ and Albanese vertex) and
$L$ for the color-profile map ($n$ rows), both acting on the space of $1$-chain
families ($n \cdot 3^{8}$ coordinates, one per colored Albanese edge, with
$|V| = 3^{8} = 6561$); membership is the solvability of
\[
M x = b, \qquad
M = \begin{pmatrix} C \\ L \end{pmatrix}, \qquad
b = (0, \dots, 0 \mid 1, \dots, 1)^{\mathsf T}.
\]
For each of the five matroids, a verified dual certificate establishes
$M(G_i) \notin \mathcal{M}_3$ on this full membership system in
$\Alb_3(M(G_i))$; the recorded ranks show the infeasibility signature
$\rank\,[M \mid b] = \rank M + 1$ in every case:

\begin{center}
\begin{tabular}{lccccc}
\hline
matroid & $(g,n)$ & rows of $M$ & columns of $M$ & $\rank M$ & $\rank\,[M \mid b]$ \\
\hline
$M(G_1)$ & $(8,16)$ & $52{,}504$ & $104{,}976$ & $48{,}274$ & $48{,}275$ \\
$M(G_2)$ & $(10,18)$ & $65{,}628$ & $118{,}098$ & $58{,}230$ & $58{,}231$ \\
$M(G_3)$ & $(9,17)$ & $59{,}066$ & $111{,}537$ & $53{,}591$ & $53{,}592$ \\
$M(G_4)$ & $(9,17)$ & $59{,}066$ & $111{,}537$ & $53{,}555$ & $53{,}556$ \\
$M(G_5)$ & $(9,17)$ & $59{,}066$ & $111{,}537$ & $53{,}612$ & $53{,}613$ \\
\hline
\end{tabular}
\end{center}

For $M(G_1)$ and $M(G_2)$ the recorded ranks exhibit the failure in its
strongest form: $\rank C = \rank M$ ($48{,}274$ and $58{,}230$ respectively,
with closedness solution spaces of dimensions $56{,}702$ and $59{,}868$). That
is, the profile rows of $L$ are linear combinations of rows of $C$, so
\emph{every} $\mathbb{Z}/3$-solution of the closedness equations has identically
zero color profile; in particular no solution with constant profile $1$ exists.
For $M(G_3), M(G_4), M(G_5)$ the rank of $C$ alone was not separately recorded,
and we state only the tabulated pair. In all five cases the dual certificate,
not the rank computation, carries the claim.

\medskip

\noindent\textbf{Minimality and pairwise independence (one argument).} For each
of the five matroids, \emph{every} single-element minor (each deletion and each
contraction) was verified to lie in $\mathcal{M}_3$ by an explicit primal
witness. For $M(G_1)$ and $M(G_2)$, whose automorphism groups act with four
orbits on the $32$ (resp.\ $36$) single-element minors, one representative per
orbit was certified (for $M(G_1)$: the four orbit classes of sizes
$12 + 12 + 4 + 4$, deletions and contractions of cube-side and apex edges); for
$M(G_3), M(G_4), M(G_5)$ all $34$ single-element minors were certified
individually. Since every proper minor of $R$ is a minor of some single-element
minor of $R$, and $\mathcal{M}_3$ is minor-closed [EGFS, Prop.~7.2], it follows
that every proper minor of each $M(G_i)$ lies in $\mathcal{M}_3$; combined with
$M(G_i) \notin \mathcal{M}_3$, each $M(G_i)$ is an excluded minor. The same fact
gives pairwise independence: a non-member of $\mathcal{M}_3$ cannot be a proper
minor of any $M(G_i)$ (all proper minors are members), so no one of the six
non-members $M(K_{3,5}), M(G_1), \dots, M(G_5)$ is a minor of another.
Distinctness as matroids also holds: rank separates all pairs except the three
rank-$9$ matroids, and their graphs are pairwise non-isomorphic and
$3$-connected (vertex-connectivity $3$ is verified with the certificates), so
their cycle matroids are pairwise non-isomorphic by Whitney's theorem. For
$M(G_1)$ the absence of a $K_{3,5}$ graph minor in $G_1$ was additionally
confirmed by direct exhaustive check.

\medskip

\noindent\textbf{A strengthening: non-membership in $\widetilde{\mathcal{M}}_3$.}
Alongside $\mathcal{M}_\ell$, [EGFS, Rmk.~8.5] introduces the larger minor-closed
class $\widetilde{\mathcal{M}}_\ell \supset \mathcal{M}_\ell$, obtained by asking
only that the color profile be nonzero at some $s \in S$ rather than constant and
nonzero, and observes that ``it is easier to show that a matroid does not lie in
$\mathcal{M}_\ell$ than it is to show that it does not lie in
$\widetilde{\mathcal{M}}_\ell$.'' For graphic matroids of biconnected graphs the
two tasks coincide, by a result of Engel and Schreieder: ``For a biconnected
graph $G$, any $\mathbb{F}_p$-solution of $\underline{R} = M(G)$ in
$\operatorname{Alb}_p(\underline{R})$ has a constant color profile
$(\lambda_s(b_s))_{s \in S} = (\lambda, \dots, \lambda)$'' [ES, Cor.~2.10],
where biconnected means that $G - v$ is connected for every vertex $v$
[ES, Prop.~2.9]. Each of $G_1, \dots, G_5$ has vertex connectivity $3$, in
particular is biconnected. Suppose some $M(G_i)$ lay in
$\widetilde{\mathcal{M}}_3$: it would admit an $\mathbb{F}_3$-solution with
profile nonzero at some $s$, hence, by the corollary, with constant profile
$(\lambda, \dots, \lambda)$, $\lambda \neq 0$; since the solution condition and
the color profile are linear in $(b_s)$, scaling by $\lambda^{-1}$ would yield a
solution with profile $\equiv 1$ --- precisely what the dual certificates above
refute. Hence $M(G_i) \notin \widetilde{\mathcal{M}}_3$ for $i = 1, \dots, 5$,
and since every proper minor lies in
$\mathcal{M}_3 \subset \widetilde{\mathcal{M}}_3$, each is an excluded minor for
$\widetilde{\mathcal{M}}_3$ as well. This strengthening is the combination of
[ES, Cor.~2.10] with the certificates of this note; it is not available from the
certificates alone. For $M(K_{3,5})$ the same conclusion
$M(K_{3,5}) \notin \widetilde{\mathcal{M}}_3$ follows from [ES, Thm.~1.7] alone
(the case $p = 3$), whose statement gives the trivial-profile form directly.

\medskip

\noindent\textbf{Verification.} The verdicts were produced by an independent
reimplementation of the membership test built solely from the definitions of
[EGFS] (Defs.~5.1, 5.3, 5.8, 8.3, 8.9), validated by reproducing every worked
value stated in [EGFS] --- the solution-space dimensions $15, 103, 35$ of
$M(K_{3,3}), M(K_5), R_{10}$; the memberships and radical distances $d = 2$ of
these three; $d = 1$ for cographic matroids; and
$M(K_{3,5}) \notin \mathcal{M}_3$ as an excluded minor. Each of the five
non-membership verdicts was reproduced by a second, independently written solver
(a distinct coordinate model for the Albanese group, a distinct matrix-assembly
order, and an independently implemented $\mathbb{F}_3$ rank computation
cross-checked against a second linear-algebra engine), with exact-integer
agreement on every recorded rank (the values tabulated above). Every verdict in
this note is certified by an explicit witness --- a dual certificate for each
non-membership, a solution vector for each minor membership --- verified by
direct arithmetic over $\mathbb{F}_3$, independently of any rank computation.

\medskip

\noindent\textbf{Code and data.} The certificate bundle accompanying this note
contains, for each of the five matroids, the dual certificate and the minor
primal witnesses of this section (for $M(G_1)$ and $M(G_2)$, one primal witness
per automorphism orbit, with the orbit accounting recorded in the bundle
manifest; for $M(G_3), M(G_4), M(G_5)$, all $34$ single-element minors
individually), the published $M(K_{3,5})$ non-membership dual as a reference
case, the graph data with integer realizations, and a single self-contained
checker, \texttt{m3\_check.py}, requiring only Python~3 with numpy. Each witness
file carries the realization matrix $A$; the checker rebuilds every Albanese
membership system from $A$ alone and verifies each witness by integer
matrix--vector products mod~$3$; the verification is independent of the search
that produced the witnesses. As a calibration, our two independent
implementations reproduce the solution dimension $56{,}381$ of the full
$K_{3,5}$ membership system at $\ell = 3$ computed in the authors' code
repository [EGFS3]. The single command \texttt{python3 m3\_check.py}
re-verifies every claim of this section; the expected final output line is
\texttt{RESULT:\ PASS} (measured: about one second and well under $1$~GB of
memory, with Python~3.13 and numpy~2.3). The checker also verifies, against
graph6 strings embedded in the checker itself, that every realization is a
reduced oriented incidence matrix of the named graph, and that the certified
minors cover all single-element deletions and contractions. The $\ell = 5$
membership witnesses of \S\ref{sec:ell5} are archived alongside, with their own
self-contained checker, in the bundle's \texttt{l5\_membership/} directory. The
bundle, together with this note, is available at
\url{https://github.com/thebookersmith/hodge-m3-excluded-minors}. Anthropic's
Claude models were used substantially in the preparation of this note and its
accompanying code; all content was reviewed by the author, who takes full
responsibility for it.

\section{Membership at $\ell = 5$}\label{sec:ell5}

The wider significance of the radical distance has recently been settled at the
level of divisibility. Engel and Schreieder [ES] prove that for a closed
subvariety $Z$ of codimension $k$ in a very general principally polarized abelian
variety of dimension $g$, the coefficient $m$ in $[Z] = m\cdot\theta_k$ (with
$\theta_k$ the minimal class) is divisible by every prime $p \le (k+1)/2$
([ES, Thm.~1.1]); equivalently $d(M(K_n))$ is divisible by all primes
$p \le (n-1)/2$ and grows at least exponentially in $n$ ([ES, Cor.~1.6]). The
prime-$p$ obstruction is realized by
$M(K_{2p+1}) \notin \mathcal{M}_p$ and $M(K_{3,2p-1}) \notin \mathcal{M}_p$
([ES, Thms.~1.5, 1.7]). For $p = 5$ the witnesses these two families provide are
$M(K_{11})$ (rank $10$, corank $45$) and $M(K_{3,9})$ (rank $11$, corank $16$): a
factor of $5$ provably appears, but only via matroids of large corank. It is then
natural to ask where, in terms of size, prime $5$ first appears. We contribute
the first data at the small-corank end.

At the time the $\ell = 5$ computation reported here was carried out, three
excluded minors for $\mathcal{M}_3$ had been explicitly identified:
$M(K_{3,5})$, $M(G_1)$, and $M(G_2)$ --- all of corank $8$.

\begin{theorem}\label{thm:ell5}
Each of $M(K_{3,5})$, $M(G_1)$, and $M(G_2)$ lies in $\mathcal{M}_5$.
\end{theorem}

Each verdict is certified by an explicit $\ell = 5$ solution --- a primal witness
$(b_s)_{s\in S}$ with $\sum_s u_s\,\partial b_s = 0$ for a basis of $\bar U$ and
constant color profile $\lambda_s(b_s) = 1$ for all $s$ --- on the full,
uncompressed membership system of [EGFS, Def.~8.3], with
$5^{8} = 390{,}625$ Albanese vertices and up to $7.03 \times 10^{6}$ colored
edges:

\begin{center}
\begin{tabular}{lccrrc}
\hline
matroid & $(g,n)$ & $|V| = 5^{8}$ & rows $(g\,5^{8}{+}n)$ & cols $(n\,5^{8})$ & verdict \\
\hline
$M(K_{3,5})$ & $(7,15)$ & $390{,}625$ & $2{,}734{,}390$ & $5{,}859{,}375$ & $\in \mathcal{M}_5$ \\
$M(G_1)$ & $(8,16)$ & $390{,}625$ & $3{,}125{,}016$ & $6{,}250{,}000$ & $\in \mathcal{M}_5$ \\
$M(G_2)$ & $(10,18)$ & $390{,}625$ & $3{,}906{,}268$ & $7{,}031{,}250$ & $\in \mathcal{M}_5$ \\
\hline
\end{tabular}
\end{center}

Every witness is checked by a self-contained verifier of the same design as the
$\ell = 3$ checker (\texttt{l5\_membership/l5\_check.py} in the archived bundle):
it rebuilds the Albanese membership system from the totally unimodular
realization alone and requires only Python with numpy.

\begin{corollary}\label{cor:dk35}
$d\bigl(M(K_{3,5})\bigr) = 6$.
\end{corollary}

\noindent The four inputs are as follows. $K_{3,5}$ contains $K_{3,3}$, hence is
non-planar, so $M(K_{3,5})$ is not cographic and
$M(K_{3,5}) \notin \mathcal{M}_2$ [EGFS, Thm.~7.1]; $M(K_{3,5}) \notin
\mathcal{M}_3$ [EGFS, Prop.~8.11]; $M(K_{3,5}) \in \mathcal{M}_5$ by
Theorem~\ref{thm:ell5}; and $M(K_{3,5}) \in \mathcal{M}_\ell$ for every prime
$\ell \ge 7$, since $\ell \ge 7 = \rank M(K_{3,5})$ and $R \in \mathcal{M}_\ell$
whenever $\ell \ge \rank R$, as observed in [EGFS, Rmk.~8.6]. Hence
$\{\ell \text{ prime} : M(K_{3,5}) \notin \mathcal{M}_\ell\} = \{2, 3\}$ and
$d(M(K_{3,5})) = \operatorname{lcm}\{2,3\} = 6$. Of the four inputs only the
$\mathcal{M}_5$-membership is new; the other three are from [EGFS].

\medskip

\noindent\textbf{Consequences at $\ell = 5$.} To our knowledge, these are the
first $\ell = 5$ membership data on the corank-$8$ slice. Three of the six
explicitly identified
excluded minors of $\mathcal{M}_3$ --- including both members of the bipartite
apex family --- do \emph{not} persist as obstructions at $\ell = 5$: each lies in
$\mathcal{M}_5$, so $5 \nmid d(R)$ for each, and in particular neither $M(G_1)$
nor $M(G_2)$ yields $30 \mid d(8)$ or $30 \mid d(10)$ (the latter divisibility
holds in any case: $M(K_{11})$ has rank $10$ and lies outside $\mathcal{M}_5$ by
[ES, Thm.~1.5], and $6 \mid d(10)$ as in \S\ref{sec:result}). Prime-$5$
obstructions do
exist in general --- by [ES], $M(K_{11})$ and $M(K_{3,9})$ lie outside
$\mathcal{M}_5$ --- so this places only these three matroids on the member side.
The $\mathcal{M}_5$-status of $M(G_3), M(G_4), M(G_5)$, the $\ell = 5$ case of
Problem~8.8, the exact corank at which a factor of $5$ first appears, and the
$\ell = 7$ status of the five new minors (whose ranks, $8$--$10$, exceed $7$)
remain open; for $M(K_{3,5})$ no prime remains undetermined
(Corollary~\ref{cor:dk35}).

\section{Remarks}\label{sec:remarks}

\noindent\textbf{Search procedure.} All five witnesses arose from a
systematic membership computation at $\ell = 3$ over a census of small non-planar
graphic matroids of corank $8$. The bipartite stratum was searched first
(motivated by the bipartite shape of the known example $K_{3,5}$) and produced
$G_1$; the same census's general stratum --- with no structural restriction ---
subsequently produced $G_3$, $G_4$, and $G_5$; $G_2$ arose from the corank-$8$
census at $11$ vertices. Within the strata searched to completion, corank
$\le 7$ contained no non-members at $\ell = 3$ at all; the concentration of all
explicitly identified excluded minors at corank $8$ reflects the census
boundary, not an
established structural fact. Likewise the search ordering reflects only
heuristics; each matroid stands on its certificate alone.

\medskip

\noindent\textbf{Outlook.} The excluded minor $M(K_{3,5})$ sits at the head
of the family $M(K_{3,n})$. Whether any of the five matroids of this note heads
an infinite family of its own is an open follow-on question, not a claim made
here; the structural diversity observed in \S\ref{sec:result} --- two bipartite
apex constructions, three sporadic non-bipartite graphs, and zero non-members
among $214$ constructed apex candidates --- suggests only that the answer to
Problem~8.8 will not be a single-family list, and the data do not yet support
more than that.

The wider context has recently shifted. The two questions of [EGFS] on the size
of the radical distance --- whether it can be arbitrarily large (Problem~8.13,
the ``second-prime''/unboundedness question) and the growth of $g \mapsto d(g)$
(Problem~8.14) --- have since been addressed by Engel and Schreieder [ES]: as
recalled in \S\ref{sec:ell5}, they prove $M(K_{2p+1}) \notin \mathcal{M}_p$ and
$M(K_{3,2p-1}) \notin \mathcal{M}_p$ for every prime $p$, whence $d(M(K_n))$ is
divisible by all primes $p \le (n-1)/2$ and grows at least exponentially. This
solves Problem~8.13 and partially answers Problem~8.14, and in particular shows
that obstructions occur at every prime: primes beyond $3$ do arise. What their
argument does not do is exhibit excluded minors --- it establishes non-membership
without minimality --- so Problem~8.8 remains open at every odd prime, and the
witnesses of [ES] are large (for $p = 5$, $M(K_{11})$ of corank $45$ and
$M(K_{3,9})$ of corank $16$). The natural refinements are therefore the
identification of excluded minors at each prime and the location of the frontier
at which each prime first appears; the $\ell = 5$ data of \S\ref{sec:ell5} is a
first probe of the small-corank end of that frontier for $p = 5$.

Finally, a methodological remark. The certificate scheme used here decides a
corank-$8$ instance at $\ell = 5$ in about eighty seconds of computation, by
exploiting the automorphism group of the underlying graph to compress the
membership system before solving. For $M(K_{3,9})$ at $\ell = 5$ (corank $16$),
where the $5^{16}$-vertex system cannot be materialized at all, a different route
succeeds: an explicit Farkas dual supported on a socle window of
$\mathbb{F}_5[t_1,\dots,t_{16}]/(t_i^5)$, giving a certificate-level verification
of the non-membership underlying [ES, Thm.~1.7] at $p = 5$. The window
certificate and its self-contained checker are included in the archived bundle
(\texttt{window\_certificates/}); the method itself will be described elsewhere.


\begin{thebibliography}{EGFS}

\bibitem[EGFS]{EGFS} P.~Engel, O.~de Gaay Fortman, S.~Schreieder,
\emph{Matroids and the integral Hodge conjecture for abelian varieties},
arXiv:2507.15704v3.

\bibitem[EGFS2]{EGFS2} P.~Engel, O.~de Gaay Fortman, S.~Schreieder,
\emph{Optimality of the Prym--Tyurin construction for $\mathcal{A}_6$},
arXiv:2512.04902 (2025).

\bibitem[EGFS3]{EGFS3} P.~Engel, O.~de Gaay Fortman, S.~Schreieder,
\emph{Solutions in Albanese graphs of regular matroids}, code repository (2025),
\url{https://github.com/philip-engel/regular_matroids}.

\bibitem[ES]{ES} P.~Engel, S.~Schreieder,
\emph{On the degree of subvarieties on abelian varieties},
arXiv:2606.31894 (2026).

\end{thebibliography}
\end{document}